\documentclass[12pt,a4paper]{amsart}

\usepackage{amssymb}
\makeatletter
\renewcommand{\@begintheorem}[2]{
\rm \trivlist \item [\hskip \labelsep {\bf #2\ \ #1.}]
                                }
\makeatother

\makeatletter

\DeclareFontFamily{U}{cyr}{}
\DeclareFontShape{U}{cyr}{m}{n}{
  <5> wncyr5 <6> wncyr6 <7> wncyr7 <8> wncyr8 <9> wncyr9 <10->
wncyr10}{}
\DeclareMathAlphabet{\mathcyr}{U}{cyr}{m}{n}

\input cyracc.def
\input epsf

\newcommand{\ZZ}{{\bf Z}}
\newcommand{\QQ}{{\bf Q}}

\newcommand{\CC}{{\bf C}}

\newcommand{\PP}{{\bf P}}

\newcommand{\lra}{\longrightarrow}

\newcommand{\cF}{{\mathcal F}}
\newcommand{\cG}{{\mathcal G}}

\newcommand{\cO}{{\mathcal O}}
\newcommand{\cP}{{\mathcal P}}

\newcommand{\tX}{{\tilde{X}}}

\newcommand{\bes}{\begin{equation*}}
\newcommand{\ees}{\end{equation*}}

\newcommand{\vP}{{P}}
\newcommand{\vQ}{{P_2}}
\newcommand{\vY}{{Y}}

\title{A remark on generalized complete intersections}
\author{Alice Garbagnati, Bert van Geemen}
\address{Dipartimento di Matematica, Universit\`a di Milano,
via Saldini 50, 20133 Milano, Italia}
\email{alice.garbagnati@unimi.it}
\address{Dipartimento di Matematica, Universit\`a di Milano,
via Saldini 50, 20133 Milano, Italia}
\email{lambertus.vangeemen@unimi.it}

\begin{document}
\begin{abstract}
We observe that an interesting method to produce non-complete intersection subvarieties,
the generalized complete intersections from L.\ Anderson and coworkers,
can be understood and made explicit by using standard Cech cohomology machinery.
We include a worked example of a generalized complete intersection Calabi-Yau threefold.
\end{abstract}
\maketitle

\section*{Introduction}
Calabi-Yau varieties, in particular those of dimension three, are of great interest in
string theory. Since there are not many general results yet on their classification,
but see \cite{Wilson}, the explicit construction of CY threefolds is a quite important enterprise.
For example, Kreuzer and Starke classified the toric fourfolds which have CY threefolds as
(anticanonical) hypersurfaces \cite{KS}, \cite{AGHJN}. Besides generalizations
to complete intersection CYs in certain ambient toric varieties,
like products of projective spaces, there are various other examples of CY threefolds
constructed with more sophisticated algebro-geometrical methods. Recent examples
include \cite{ito}, \cite{CGKK}, \cite{kapustka}.

In the recent paper \cite{anderson}, L.\ Anderson, F.\ Apruzzi, X.\ Gao, J.\ Gray and S-J.\ Lee
found a very nice method to construct many more CY threefolds. The basic idea is to take
a hypersurface $Y$ in an ambient variety $P$ and to consider hypersurfaces $X$ in $Y$.
These hypersurfaces need not be complete intersections in $P$, that is, there need not exist
two sections of two line bundles on $P$ whose common zero locus is $X$.
There are various generalizations of this method, but we will stick to this
basis case.
As in \cite{anderson}, we refer to these varieties as generalized complete intersections
(gCIs).

A particularly interesting and accessible case that was found and studied by
Anderson and coworkers is when the ambient variety is a product of two varieties,
one of which is $\PP^1$, so $\vP=\vQ\times\PP^1$.
The variety $P_2$ they consider is a product of projective spaces,
but this is not essential, one could consider any toric variety or even more general cases.
The factor $\PP^1$ is important since there are line bundles on $\PP^1$ with non-trivial
first cohomology group and this is essential to find generalized complete intersections.
We review this construction in Section \ref{genset}.

We provide a proposition, proven with standard Cech cohomology methods,
that allows one, under a certain hypothesis, to find three equations
(more precisely, three sections of three line bundles on
$P$) that define $X$. In Section \ref{exa} we work out a detailed example, with explicit
equations, of a CY threefold which was already considered in \cite{anderson}.
The explicit example $X$ has an automorphism of order two and the quotient
of $X$ by the involution provides, after desingularization, another CY threefold.
More generally, we think that among the gCIs found in \cite{anderson} one could find
more examples of CY threefolds with non-trivial automorphisms. It might be hard though to
implement a systematic search as was done in \cite{CGL} for complete intersection CY
threefolds in products of projective spaces.
We did not find new CY threefolds with small Hodge numbers
(see \cite{CCM} for an update on these),
but the gCICY seem to be a promising class of CYs to search for these.
The recent paper \cite{berglund} by Berglund and H\"ubsch provides further techniques to deal
with gCICYs whereas \cite{anderson2} explores string theoretical aspects of gCICYs.

\section{The construction of generalized complete intersections}

\subsection{The general setting}\label{genset}
Let $\vQ$ be a projective variety of dimension $n$ and let $\vP:=\vQ\times\PP^1$.
We denote the projections to the factors of $\vP$ by $\pi_1,\pi_2$ respectively.
For  a coherent sheaf $\cF$ on $\vQ$ and an integer $d$
we define a coherent sheaf on $\vP$ by:
$$
\cF[d]\,:=\,\pi_1^*\cF\otimes \pi_2^*\cO_{\PP^1}(d)~.
$$
The K\"unneth formula gives
$$
H^r(\vP,\cF[d])\,=\,\oplus_{p+q=r}\, H^p(\vQ,\cF)\otimes H^q(\PP^1,\cO_{\PP^1}(d))~.
$$
Recall that the only non-zero cohomology of $\cO_{\PP^1}(d)$ is:
$h^0(\cO_{\PP^1}(d))=h^1(\cO_{\PP^1}(-2-d))=d+1$ for $d\geq 0$ and a basis for
$H^0(\cO_{\PP^1}(d))$ is given by the monomials $z_0^iz_1^{d-i}$, $i=0,\ldots,d$,
where $(z_0:z_1)$ are the homogeneous coordinates on $\PP^1$.

Let $L$ be a line bundle on $\vQ$ and assume that $L[d]$,
for some $d\geq 1$, has a non-trivial global section $F$.
Using the K\"unneth formula,
we can write $F=\sum_if_iz_0^iz_1^{d-i}$ for certain sections $f_i\in H^0(\vQ,L)$.
Let $\vY=(F)$ be the zero locus of $F$ in $\vP$.
We assume that $\vY$ is a (reduced, irreducible) variety,
although this will not be essential in this section.

To define a codimension two subvariety of $\vP$, we consider another line bundle $M$ on
$\vQ$. The K\"unneth formula shows that $M[-e]$ has no global sections if $e\geq 1$.
But upon restricting to $Y$, the vector space $H^0(\vY,M[-e]_{|\vY})$
could still be non-trivial.
In fact, from the exact sequence
\begin{equation}\label{koszul}
0\,\lra\,(L^{-1}\otimes M)[-d-e]\,\stackrel{F}{\lra}M[-e]\,\lra\,M[-e]_{|\vY}\,\lra\,0
\end{equation}
we deduce the exact sequence
\begin{equation}\label{F1}
0\,\lra H^0(\vY,M[-e]_{|\vY})\,\stackrel{d^0}{\lra}\,H^1(P,(L^{-1}\otimes M)[-d-e])\,
\stackrel{F_1}{\lra} H^1(P,M[-e])
\end{equation}
thus $H^0(Y,M[-e]_{|\vY})\cong \ker(F_1)$,
where we denote by $F_1$ the map induced by multiplication by
$F$ on the first cohomology groups.
Since now
$$
H^1(P,(L^{-1}\otimes M)[-d-e])\,\cong\,H^0(P_2,L^{-1}\otimes M)\otimes
H^1(\PP^1,\cO_{\PP^1}(-d-e))\qquad (d+e\geq 2)
$$
the domain of $F_1$ is non trivial if and only if
$h^0(P_2,L^{-1}\otimes M)\neq 0$.
So for suitable choices of line bundles on $P_2$
we might find interesting, non-complete intersection, codimension two
subvarieties of $\vP$ in this way. In the proof of Proposition \ref{prop}
we explain how to compute $F_1$.

\subsection{Example}\label{exaPn}
Let $\vQ=\PP^n$, $L=\cO_{\PP^n}(k)$, $M=\cO_{\PP^n}(k+l)$ with $l\geq0$,
let $d\geq 1$ and $e=1$. Then $h^0(\PP^n,L^{-1}\otimes M)=h^0(\cO_{\PP^n}(l))\neq 0$
so $h^1(P,(L^{-1}\otimes M)[-d-e])\neq 0$,
but $h^1(P,M[-e])=0$ since $\cO_{\PP^1}(-1)$ has no cohomology. Thus
$H^0(Y,M[-e]_{|\vY})\cong H^1(P,(L^{-1}\otimes M)[-d-e])$ is indeed non-trivial.

\subsection{Generalized complete intersections}\label{exeq}
Given a variety $\vY\subset P$ that is the zero locus of $F\in H^0(L[d])$
as in Section \ref{genset}, and given a global section
$\tau\in H^0(M[-e]_{|\vY})$,
its zero locus $X:=(\tau)\subset\vY$ is called a generalized complete intersection.

The scheme $X$ may not be defined
by two global sections $\sigma_1,\sigma_2$ of line bundles $L_1$, $L_2$ on $\vP$.
However in certain cases we can find three sections of line bundles on $P$ which define $X$:

\subsection{Proposition}\label{prop}
Let $F\in H^0(\vP,L[d])$, let $\vY=(F)$, let
$\tau\in H^0(\vY,M[-e]_{|\vY})$ with $d,e\geq 1$ be as above and assume that
$H^1(\vQ,L^{-1}\otimes M)=0$.

Then there are two global sections $G,H\in H^0(\vP,M[d-1])$ such that the
generalized complete intersection subscheme $X$ of
$\vP$ defined by $\tau$ in $\vY$ can also be defined as
$$
X\,=\,\{x\in\vP:\,F(x)\,=\,G(x)\,=\,H(x)\,=\,0\}
$$
(the equality is of schemes). Moreover, there is a global section
$A\in H^0(\vP,(L^{-1}\otimes M)[d+e-2])$
such that
$AF=z_1^{d+e-1}G+z_0^{d+e-1}H$, so that on the open subset of $\vQ\times \PP^1$ where
$z_0\neq 0$ the subscheme $X$ of $\vP$ is defined by the two equations $F=G=0$.

\

\noindent {\bf Proof.}
We use Cech cohomology to make the isomorphism $H^0(\vY,M[-e]_{|\vY})\cong \ker(F_1)$,
see exact sequence (\ref{F1}),
explicit. Let $U_i\subset\PP^1$ be the open subset where $z_i\neq 0$.
For a coherent sheaf $\cG$ on $\PP^1$ we have the exact sequence
$$
0\,\lra\,H^0(\PP^1,\cG)\,\lra\,\cG(U_0)\,\oplus\,\cG(U_1)\,\stackrel{\delta}{\lra}\,
\cG(U_0\cap U_1)\,\lra \,H^1(\PP^1,\cG)\,\lra 0~,
$$
where $\delta(t_0,t_1)=t_0-t_1$.
The cohomology groups we consider are computed with the K\"unneth formula.
Note that after tensoring this exact sequence by a vector space $W$,
we obtain that $W\otimes H^0(\PP^1,\cG)=\ker(1_W\otimes\delta)$ and
$W\otimes H^1(\PP^1,\cG)=\operatorname{coker}(1_W\otimes\delta)$.

For an affine open subset $V\subset\PP^1$,
the cohomology of the exact sequence (\ref{koszul}) on $\vQ\times V$ gives the exact
sequence, where we extend $M[-e]_{|\vY}$ by zero to $\vQ\times V$,
$$
H^0(\vQ\times V,M[-e])\,\lra\,H^0(\vQ\times V,M[-e]_{|\vY})\,\lra\,
H^1(\vQ\times V,(L^{-1}\otimes M)[-d-e])~.
$$
The K\"unneth formula, combined with the assumption $H^1(P_2,L^{-1}\otimes M)=0$
and the fact that $H^1(V,\cF)=0$ for any coherent sheaf $\cF$ since $V$ is affine,
implies that the last group is zero.

Taking $V=U_0,U_1$, the exact sequence
(\ref{koszul}) on $\vQ\times V$ thus gives two exact sequences whose
 sum (term by term) is
{\renewcommand{\arraystretch}{1.3}
\begin{align}\label{sequence 1}\begin{array}{cccccc}
0\lra\oplus_{i=0}^1 H^0(L^{-1}\otimes M)\otimes (\cO_{\PP^1}(-d-e)(U_i))
\stackrel{F}{\lra}\\
\oplus_{i=0}^1
H^0(M)\otimes(\cO_{\PP^1}(-e)(U_i))\lra \oplus_{i=0}^1
(M[-e]_{|\vY})(\vQ\times U_i)\lra0~.
\end{array}\end{align}
Similarly taking $V=U_0\cap U_1$ one has the exact sequence:
\begin{align}\label{sequence 2}\begin{array}{cccccc}
0\lra H^0(L^{-1}\otimes M)\otimes (\cO_{\PP^1}(-d-e)(U_0\cap U_1))
\stackrel{F}{\lra}\\
 H^0(M)\otimes (\cO_{\PP^1}(-e)(U_0\cap U_1)) \lra
(M[-e]_{|\vY})(\vQ\times( U_0\cap U_1))\lra0~.
\end{array}
\end{align}
}
Next we use the Cech boundary map $\delta$ to map
sequence \eqref{sequence 1} to sequence \eqref{sequence 2} and we obtain
a commutative diagram
with three complexes as columns.
The first two columns are Cech complexes for the covering $\{U_i\}_{i=0,1}$ of $\PP^1$,
their cohomology groups are respectively
{\renewcommand{\arraystretch}{1.3}
$$
\begin{array}{ccc}
H^0(L^{-1}\otimes M)\otimes H^q(\cO_{\PP^1}(-d-e))
&\cong&
H^q(P,(L^{-1}\otimes M)[-d-e]),\\
\qquad
H^0(M)\otimes H^q(\cO_{\PP^1}(-e))
&\cong&
H^q(P,M[-e]), \quad(q=0,1)~.
\end{array}
$$
}
The zero-th cohomology group of the last column is $H^0(Y,M[-e]_{|Y})$.
So we conclude that the maps $q$ and $F_1$ can be computed with the long exact cohomology
sequence associated to this diagram.

We observe, but will not use, that the K\"unneth formula implies that
$H^2(P,(L^{-1}\otimes M)[-d-e])=0$ and thus the cohomology sequence
of (\ref{koszul}) gives a six term exact sequence with the zero-th and first cohomology groups.
The first 5 terms are the same as those of the long exact sequence associated to the
diagram, so we conclude that the first cohomology group of the last column
is $H^1(Y,M[-e]_{|Y})$.

Given $\tau\in H^0(\vY,M[-e]_{|\vY})$, let $q:= d^0(\tau)\in \ker(F_1)$.
Since the first row (\ref{sequence 1}) of the complex is exact,
the section $\tau$ is locally given by restricting sections
$\tau_i\in M[-e](\vQ\times U_i)$ to $\vY$.
By the snake lemma, they satisfy $\tau_0-\tau_1=Fq$ on $\vQ\times (U_0\cap U_1)$,
in particular $\tau_0=\tau_1$ on $\vY\cap (\vQ\times (U_0\cap U_1))$ since $F=0$ on $Y$.

The images of the $z_0^{-j}z_1^{-d-e+j}\in \cO_{\PP^1}(-d-e)(U_0\cap U_1)$, $j=1,\ldots,d+e-1$,
form a basis of $H^1(\PP^1,\cO(-d-e))$. A cohomology class
$q\in H^1(\vP,(L^{-1}\otimes M)[-d-e])
\cong H^0(\vQ,L^{-1}\otimes M)\otimes H^1(\PP^1,\cO(-d-e))$
can thus be represented by $q=\sum_j q_jz_0^{-j}z_1^{-d-e+j}$ with
$q_i\in H^0(\vQ,L^{-1}\otimes M)$.
Let $F=\sum_if_iz_0^iz_1^{d-i}$, where $f_i\in H^0(P_2,L)$,
then $Fq$ is homogeneous of degree $d-(d+e)=-e$ and
it is a sum of terms $r_kz_0^kz_1^{-e-k}$ with $r_k\in H^0(P_2,M)$.
Writing
$$
Fq=\,\sum_{k=-d-e+1}^{d-1} r_kz_0^kz_1^{-e-k}\,=
\left(\sum_{k=-d-e+1}^{-e} r_kz_0^kz_1^{-e-k} \right)+
\left(\sum_{k=-e+1}^{-1} r_kz_0^kz_1^{-e-k} \right)+
\left(\sum_{k=0}^{d-1}r_kz_0^kz_1^{-e-k} \right),
$$
the first summand lies in $M[-e](\vQ\times U_0)$ (where $z_0\neq 0$)
and the last summand
lies in $M[-e](\vQ\times U_1)$, we denote these summand by $\tau_0$ and $-\tau_1$ respectively.
The middle summand has monomials $z_0^az_1^b$ with both $a,b<0$.
Thus $Fq$ represents a class in $q'\in H^1(\vP,M[-e])$, which is the same as the class
represented by the middle summand. By definition, one has $q'=F_1(q)$ and thus $q\in\ker(F_1)$
when all coefficients $r_k$, $k=-e+1,\ldots,-1$, are zero.

Since $q\in \ker(F_1)$ this middle summand is zero, so that
$Fq=\tau_0-\tau_1$ as desired.
Now we define $G:=z_0^{d+e-1}\tau_0$ and $H:=-z_1^{d+e-1}\tau_1$
so that all their monomials $z_0^az_1^b$ have $a,b\geq 0$ and $a+b=d-1$,
thus both $G,H\in H^0(\vP,M[d-1])$.
Then $(z_0z_1)^{d+e-1}Fq=z_1^{d+e-1}G+z_0^{d+e-1}H$
and with $A:=(z_0z_1)^{d+e-1}q\in H^0(\vP,(L^{-1}\otimes M)[d+e-2])$
we find the desired relation.
\qed

\subsection{Example}\label{exa2}
With the choices of $P_2,L,M$ as in Example \ref{exaPn},
and if $X$ is a smooth variety (of dimension $n-1$), then $H^1(\vQ,L^{-1}\otimes M)=
H^1(\PP^n,\cO_{\PP^n}(l))=0$, for any $l$, if $n>1$.
The adjunction formula implies that
$X$ has trivial canonical bundle if we choose $l=n+1-2k$ and $d=3$.
In that case $P=\PP^n\times\PP^1$ and $F$ is homogeneous of bidegree $(k,3)$
whereas $G,H$ have bidegree $(n+1-k,2)$.

\subsection{A fibration on $X$}\label{fib1}
Given $X$ as in the proposition, the projection $\pi_2: P_2\times\PP^1 \rightarrow\PP^1$
restricts to $X$ to give a fibration denoted by $\pi_2:X\rightarrow\PP^1$.
For a point $p=(z_0:z_1)\in\PP^1$, we denote by
$F_p\in H^0(P_2,L)$, $H_p\in H^0(P_2,M)$
the restrictions of $F$ and $H$ to the fiber $X_p$.
The equation
$AF=z_1^{d+e-1}G+z_0^{d+e-1}H$
shows that if $z_1\neq 0$ then $F_p$ and $H_p$
define the fiber $X_p$, which is thus a complete intersection in $P_2$.

\subsection{Example}
This example illustrates that $X$, as in Proposition \ref{prop},
might be reducible, even if $h^0(Y,M[-e]_{|Y})$ is rather large. The example is taken
from \cite[Table 4]{anderson}, third item (with $i=2$) where it is in fact observed
that no smooth varieties arise in that case.
We take
$$
\vQ:=\PP^2\times \PP^1\times\PP^1~,\qquad L:=\cO(0,1,1),\quad M:=\cO(3,1,1)~,
\qquad d=4,\;e=-2~.
$$
Notice that $H^1(\vQ,L^{-1}\otimes M)=H^1(\PP^2\times \PP^1\times\PP^1,\cO(3,0,0))=0$ by the
K\"unneth formula, so we can, but will not, apply Proposition \ref{prop}.
Since $h^1((L^{-1}\otimes M)[-d-e])=h^1(\cO(3,0,0)[-6])=10\cdot1\cdot1\cdot5=50$
and $h^1(M[-e])=10\cdot2\cdot2\cdot1=40$, we find $h^0(M[-e]_{|Y})\geq 10$.
We will show that, for general $Y$, $h^0(M[-e]_{|Y})=10$ but that all sections of
$M[-e]_{|Y}$ define reducible subvarieties of $Y$.

Due to the first zero in $L=\cO(0,1,1)$, the variety $Y$ is a product,
$Y=\PP^2\times S\subset P$,
with $S\subset (\PP^1)^3$ the surface defined by a section of $\cO(1,1,4)$.
Then we have
$h^0(M[-e]_{|Y})=h^0(\PP^2\times S,\pi_1^*\cO_{\PP^2}(3)\otimes\pi_2^*\cO_S(1,1,-2))$
and using the K\"unneth formula we find
$h^0(M[-e]_{|Y})=h^0(\cO_{\PP^2}(3))h^0(\cO_S(1,1,-2))=10h^0(\cO_S(1,1,-2)$.
The exact sequence
$$
0\,\lra \cO_{(\PP^1)^3}(0,0,-6)\,\stackrel{f}{\lra}\,\cO_{(\PP^1)^3}(1,1,-2)\,
\stackrel{}{\lra}\,\cO_S(1,1,-2)\,\lra\,0~,
$$
where $f$ is the equation of $S$, shows that (with $f_1$ the map induced by $f$ on $H^1$):
$$
h^0(\cO_S(1,1,-2))=\dim
\ker\Big(f_1:H^1(\cO_{(\PP^1)^3}(0,0,-6))\rightarrow H^1(\cO_{(\PP^1)^3}(1,1,-2))\Big)~.
$$
Since these spaces have dimensions $1\cdot1\cdot 5=5$ and $2\cdot 2\cdot 1=4$ respectively,
one expects $h^0(\cO_S(1,1,-2))=1$.
In that case any section $\tau\in H^0(M[-e]_{|Y})$
would be the product $\tau=gs$ with $g\in H^0(\cO_{\PP^2}(3))$
and $s\in H^0( \cO_S(1,1,-2))$ the unique (up to scalar multiple) section, hence $X$ would
be reducible.

To see that indeed $h^0(\cO_S(1,1,-2))=1$ for a general equation $f$,
take a smooth (genus one) curve $C$ of
bidegree $(2,2)$ in $\PP^1\times\PP^1$ and choose eight distinct points on $C$
which are not cut out by another curve of bidegree $(2,2)$.
As curves of bidegree $(1,4)$ depend on $2\cdot 5=10$
parameters, we can find two polynomials $g_0,g_1$
of bidegree $(1,4)$ such that $g_0=g_1=0$ consists of these eight points on $C$.
Take $f=x_0g_0+x_1g_1$ with $(x_0:x_1)\in\PP^1$, the first copy of $\PP^1$ in $(\PP^1)^3$,
and the $g_i$ on the last two copies of $\PP^1$.
The surface $S\subset (\PP^1)^3$ defined by $f$ is thus the blow up of
$\PP^1\times\PP^1$ in the eight points where $g_0=g_1=0$.
The adjunction formula shows that the
line bundle $\cO_S(1,1,-2)$ is the anticanonical bundle of $S$.
The effective anticanonical divisors are the strict transforms of bidegree $(2,2)$-curves
on passing through these eight points. Hence
the strict transform of $C$ in $S$ will be the unique effective anticanonical divisor
on $S$ and therefore $h^0(\cO_S(1,1,-2))=1$.

\section{An example: a generalized complete intersection Calabi-Yau threefold}\label{exa}

\subsection{} We illustrate the use of Proposition \ref{prop} (and its proof) for the
generalized complete intersection Calabi Yau discussed in \cite[Section 2.2.2]{anderson}.
We also consider an explicit example which has a non-trivial involution
and we compute the
Hodge numbers of a desingularization of the quotient threefold which is again a CY.

\subsection{The varieties $\vQ$ and $\vY$}\label{exvar}
We consider the case that $\vQ=\PP^4$, we choose the line bundle
$L:=\cO_{\PP^4}(2)$ and we let $d=3$.
Then the line bundle $L[d]=\cO_{\vP}(2,3)$ is very ample on $\vP=\PP^4\times\PP^1$ and thus
a general section $F$ will define a smooth fourfold $\vY$ of $\vP$.
To obtain a CY threefold in $\vY$, we consider global sections of the anticanonical bundle
of $\vY$. By adjunction, $\omega_Y=(\cO_\vP(-5,-2)\otimes\cO_\vP(2,3))_{\vY}=\cO_\vY(-3,1)$.
Thus we take $M=\cO_{\PP^4}(3)$ and $e=1$, so that $M[-e]_{|Y}=\cO_\vY(3,-1)=\omega_\vY^{-1}$.
As the $H^1$ of any line bundle on $\PP^4$ is trivial, we can use (the proof of)
Proposition \ref{prop} to find polynomials $G,H\in H^0(\vP,\cO_\vP(3,2))$ which together with
$F$ define a generalized complete intersection $X$.

As in Example \ref{exaPn}, we get
$$
H^0(\cO_\vY(3,-1))\,\stackrel{\cong}{\lra}\,H^1(\cO_\vP(1,-4))~.
$$
To find explicit elements of $H^0(\cO_\vY(3,-1))$, we write the defining equation of $Y$ as
$$
F\,=\, P_0z_0^3\,+\,P_1z_0^2z_1\,+\,P_2z_0z_1^2\,+\,P_3z_1^3\qquad(\in H^0(P,\cO_{\vP}(2,3)))~,
$$
with $P_i\in H^0(\PP^4,\cO(2))$ homogeneous polynomials of degree two in $y=(y_0:\ldots:y_4)$.
As $H^1(\cO_\vP(1,-4))\cong H^0(\cO_{\PP^4}(1))\otimes H^1(\cO_{\PP^1}(-4))$, a basis of this
$5\cdot 3=15$ dimensional vector space are the products of one of
$y_0,\ldots,y_4$ with one of $z_0^{-3}z_1^{-1}, z_0^{-2}z_1^{-2},z_0^{-1}z_1^{-3}$.
Thus any class $q\in H^1(\cO_\vP(1,-4))$ has a representative
$$
q\,=\,Q_0z_0^{-3}z_1^{-1}\,+\,Q_1z_0^{-2}z_1^{-2}\,+\,Q_2z_0^{-1}z_1^{-3}
\qquad(\in H^1(\cO_\vP(1,-4)))~,
$$
with linear forms $Q_i\in H^0(\PP^4,\cO(1))$.
As in the proof of Proposition \ref{prop} we must write:
$$
Fq\,=\,\tau_0-\tau_1,\qquad G\,:=\,z_0^3\tau_0,\qquad H\,:=\,-z_1^3\tau_1~,
$$
with $\tau_i\in \cO_P(3,-1)(\PP^4\times U_i)$.
So we find
{\renewcommand{\arraystretch}{1.3}
$$
\begin{array}{rcl}
G&=&\phantom{-(\,}z_0^2(P_1Q_0+P_2Q_1+P_3Q_2)\,+\,z_0z_1(P_2Q_0+P_3Q_1)\,+\,z_1^2P_3Q_0~,
\\
H&=&-\Bigl(\,z_0^2P_0Q_2\,+\,z_0z_1(P_0Q_1\,+\,P_1Q_2)\,+\,z_1^2(P_0Q_0+P_1Q_1+P_2Q_2)
\,\Bigr)~.
\end{array}
$$
}

\subsection{The base locus of $|-K_Y|$}
In Section \ref{exvar} we showed how to find the global sections of $\omega^{-1}_Y=\cO_Y(3,-1)$
explicitly, locally such a section is given by the polynomials $G$ and $H$.
From the formula for $F$ we see that if
$x\in \PP^4$ and $P_0(x)=\ldots=P_3(x)=0$, then the curve $\{x\}\times\PP^1$ lies in
$Y$. This curve also lies in the zero loci of $G$ and $H$, for any choice of $Q_0,Q_1,Q_2
\in H^0(\cO_{\PP^4}(1))$, hence it lies in the base locus of anticanonical system
$|-K_Y|$. Since the four quadrics $P_i=0$ in $\PP^4$ intersect in at least $2^4$ points,
counted with multiplicity, we see that this base locus is non-empty.
Thus we cannot use Bertini's theorem
to guarantee that there are smooth CY threefolds $X\subset Y$, but we resort to an
explicit example, see below.

\subsection{The CY threefold $X$}\label{exX}
To obtain an explicit example, we choose
$$
P_0\,:=\,y_0^2+y_1^2+y_2^2+y_3^2+y_4^2,
\quad
P_1\,:=\,y_0^2+y_4^2,
\quad
P_2\,:=\,y_1^2+y_3^2,
\quad
P_3\,:=\,y_0^2+y_1^2-y_2^2-y_3^2-y_4^2,
$$
and
$$
Q_0:=y_0,\qquad  Q_1:=y_1,\qquad Q_2:=y_2~.
$$
Using a computer algebra system (we used Magma \cite{magma}), one can verify that
$Y:=(F=0)$ and
$X:=(F=G=H=0)$ are smooth varieties in $\vP$. The variety $X$ is a Calabi-Yau threefold
since it is an anticanonical divisor on $Y$. In \cite[(2.27), (2.28)]{anderson}
one finds that the Hodge numbers of $X$ are $(h^{1,1}(X),h^{2,1}(X))=(2,46)$,
in particular, $h^2(X)=2$, $h^3(X)=94$.

\subsection{Parameters}
The CY threefold $X$ is defined by a section $F\in H^0(\vP,\cO_\vP(2,3))$
and a section $\tau\in H^0(\vY,\cO_\vY(3,-1))$. The first is a vector space of dimension
$$
h^0(\vP,\cO_\vP(2,3))\,=\,h^0(\PP^4,\cO_{\PP^4}(2))\cdot h^0(\PP^1,\cO_{\PP^1}(3))
\,=\,15\cdot 4=60~,
$$
whereas the second has dimension $15$. The group $GL(5,\CC)\times GL(2,\CC)$ acts on
$H^0(\cO_\vP(2,3))$
and has dimension $5^2+2^2=29$. The subgroup
of elements $(\lambda I_5,\mu I_2)$ with $\lambda^2\mu^3=1$ acts trivially, so we get
$60-28=32$ parameters for $\vP$ and next $15-1=14$ parameters for $\tau$,
so we do get $32+14=46=h^{2,1}(X)$ parameters for $X$.  So the general
deformation of $X$ seems to be  again a gCICY of the same type as $X$. (In
\cite{anderson}, just below (2.28), the dependence of $X$ on $\cP$, which gives $32$
parameters, seems to have been overlooked.)

\subsection{A CY quotient}\label{exquo}
A well-known method to obtain Calabi-Yau threefolds is to consider desingularizations
of quotients of such threefolds by finite groups, see for example \cite{CGL}.
In the example above, we see that $X\subset\PP^4\times\PP^1$ has a
subgroup $(\ZZ/2\ZZ)^2\subset Aut(X)$ given by the sign changes of $y_3$ and $y_4$.
We consider the involution
$$
\iota:\,X\,\longrightarrow\,X,\qquad
\Big((y_0:\ldots:y_4),(z_0:z_1)\Big)\longmapsto \Big((y_0:y_1:y_2:-y_3:-y_4),(z_0:z_1)\Big)~.
$$
Its fixed point locus has two components,
one defined by $y_3=y_4=0$ and the other by $y_0=y_1=y_2=0$ in  $X$.
The first is a curve in $\PP^2\times\PP^1\subset \vP$, which is smooth, irreducible and reduced
of genus $8$  according to Magma.
Similarly, the other component is a genus $2$ curve in $\PP^1_{(y_3:y_4)}\times\PP^1_{(z_0:z_1)}
\subset \vP$. In fact, only $F=0$ provides a non-trivial equation for this curve since
$y_0=y_1=y_2=0$ implies $Q_0=Q_1=Q_2=0$ and hence $G=H=0$ on this $\PP^1\times\PP^1$.
As $F=0$ defines a smooth curve of bidegree $(2,3)$ in $\PP^1\times\PP^1$,
this curve has genus $(2-1)(3-1)=2$.

In particular, the singular locus of the quotient $X/\iota$
consists of two curves of $A_1$-singularities.
Since the fixed point locus $X^\iota$ consists of two curves, we conclude that
locally on $X$ the involution is given by $(t_1,t_2,t_3)\mapsto(-t_1,-t_2,t_3)$
in suitable coordinates.
Hence $\iota$ acts trivially on the nowhere vanishing holomorphic 3-form on
the CY threefold $X$.
Thus the blow up $Z$ of $X/\iota$ in the singular locus will again be a CY threefold.

We determine the Hodge numbers of $Z$.
To do so, it is more convenient to consider the blow up $\tX$
of $X$ in the fixed point locus $X^\iota$. The involution extends to an involution
$\tilde{\iota}$ on $\tilde{X}$, the fixed point set of $\tilde{\iota}$ consists
of the two exceptional divisors and the quotient $\tX/\tilde{\iota}$ is the same $Z$.
Moreover, $H^i(Z,\QQ)\cong H^i(\tX,\QQ)^{\tilde{\iota}}$,
the $\tilde{\iota}$-invariant subspace.

Standard results on the blow up of smooth varieties in smooth subvarieties
(cf. \cite[Thm 7.31]{Voisin})
show that
$h^2(\tX)=h^2(X)+2=4$ (due to the two exceptional divisors over the two fixed curves) and
$h^3(\tX)=h^3(X)+2\cdot8 +2\cdot 2=114$
(the contribution of the $H^1$ of the fixed curves to $H^3$ of the blow up).
The Lefschetz fixed point formula for $\tilde{\iota}$ gives
$$
\chi(\tX^{\tilde{\iota}})\,=\,\sum_{i=0}^6(-1)^itr(\tilde{\iota}^*|H^i(\tX,\QQ))~.
$$
Notice that $\tilde{\iota}^*$ is the identity on $H^0,H^2,H^4,H^6$, in particular
$h^2(Z)=\dim H^2(\tX,\QQ)^{\tilde{\iota}}=4$.
The fixed points of $\tilde{\iota}$ are the two exceptional divisors,
these are $\PP^1$-bundles over the exceptional curves hence
$$
2(2-2\cdot 2)+2(2-2\cdot 8)\,=\,1-0+4-t_3+4-0+1\quad\Longrightarrow\quad t_3=42~.
$$
If the $+,-$ eigenspaces of $\tilde{\iota}$
on $H^3(\tX,\QQ)$ have dimensions $a,b$ respectively, then $a+b=114$ and $a-b=42$, thus
$a=78$ and $a=\dim H^3(\tX,\QQ)^{\tilde{\iota}}=h^3(Z)$.
As $Z$ is a CY threefold it has $h^{3,0}(Z)=1$ and thus $h^{2,1}(Z)=(78-2)/2=38$.
Other examples of CY threefolds with $(h^{1,1},h^{2,1})=(4,38)$ are already known.

\subsection{A (singular) projective model of $Z$}
The fibers of $\pi_2:X\rightarrow\PP^1$ are K3 surfaces, complete intersections of a quadric
and a cubic hypersurface in $\PP^4$. The involution $\iota$ on $X$ restricts
to a Nikulin involution on each smooth fiber. The quotient of such a fiber by the involution
will in general be isomorphic to a K3 surface in $\PP^2\times\PP^1$, defined by an equation
of bidegree $(3,2)$ (see \cite[Section 3.3]{GS}). Using the same method as in
that reference, we found that the rational map
$$
\PP^4\times\PP^1\,--\rightarrow\,\PP^2\times\PP^1\times\PP^1,\quad
\Big((y_0:\ldots:y_4),(z_0:z_1)\Big)\,\longmapsto\,
\Big((y_0:y_1:y_2),(y_3:y_4),(z_0:z_1)\Big)
$$
factors over $X/\iota$ and the image, defined by an equation of multidegree $(3,2,2)$,
is birational with $Z$. Using the explicit equation for the image and Magma,
we found that the image has $38$ singular points.


\begin{thebibliography}{AAGGL}
\bibitem[AAGGL]{anderson} L.B.\ Anderson,  F.\ Apruzzi, X.\ Gao, J.\ Gray, S-J.\ Lee,
{\it A New Construction of Calabi-Yau Manifolds: Generalized CICYs},
(arXiv:1507.03235),  Nuclear Phys.\ B {\bf 906} (2016) 441--496.

\bibitem[AAGGL2]{anderson2}  L.B.\ Anderson,  F.\ Apruzzi, X.\ Gao, J.\ Gray, S-J.\ Lee,
{\it 
Instanton Superpotentials, Calabi-Yau Geometry, and Fibrations},
(arXiv:1511.05188), Phys.\ Rev.\ D {\bf 93} (2016), no.\ 8, 086001.

\bibitem[AGHJN]{AGHJN}
R.\ Altman, J.\ Gray, Y-H.\ He,  V.\ Jejjalaf, B.D.\ Nelson,
{\it A Calabi-Yau database: threefolds constructed from the Kreuzer-Skarke list},
J.\ High Energy Phys.\ 2015, no.\ 2, 158.

\bibitem[BH]{berglund} P.\ Berglund, T.\ H\"ubsch,
{\it On Calabi-Yau generalized complete intersections from Hirzebruch varieties
and novel K3-fibrations},
arXiv:1606.07420.

\bibitem[CCM]{CCM} P.\ Candelas, A.\ Constantin, and C.\ Mishra,
{\it Calabi-Yau Threefolds With Small Hodge Numbers},
arXiv:1602.06303.

\bibitem[CGL]{CGL} A.\ Constatin, J.\ Gray, A.\ Lukas,
{\it Hodge numbers for all CICY Quotients}, arXiv:1607.01830.

\bibitem[CGKK]{CGKK} S.\ Coughlan, L.\ Golebiowski, G.\ Kapustka, M.\ Kapustka,
{\it Arithmetically Gorenstein Calabi-Yau threefolds in $P^7$},
Electron.\ Res.\ Announc.\ Math.\ Sci.\ {\bf 23} (2016) 52--68.


\bibitem[GS]{GS} B.\ van Geemen, A.\ Sarti,
{\it Nikulin involutions on K3 surfaces},
Math.\ Z.\ {\bf 255} (2007) 731--753.

\bibitem[IMOU]{ito} A.\ Ito, M.\ Miura, S.\ Okawa, K.\ Ueda,
{\it Calabi-Yau complete intersections in $G_2$-Grassmannians},
arXiv:1606.04076.


\bibitem[KK]{kapustka} G.\ Kapustka, M. Kapustka,
{\it Calabi-Yau threefolds in $P^6$},
Ann.\ Mat.\ Pura Appl.\ {\bf 195} (2016) 529--556.


\bibitem[KS]{KS}  M.\ Kreuzer and H.\ Skarke,
{\it Complete classification of reflexive polyhedra in four-dimensions},
Adv.\ Theor.\ Math.\ Phys.\ {\bf 4}, 1209 (2002).


\bibitem[M]{magma} W.\ Bosma, J.\ Cannon, and C.\ Playoust,
{\it The Magma algebra system. I. The user language}, J.\ Symbolic
Comput.\ {\bf 24} (1997) 235--265.

\bibitem[V]{Voisin}
C.\ Voisin,
Hodge Theory and Complex Algebraic Geometry I. Cambridge University Press 2002.

\bibitem[W]{Wilson} P.M.H.\ Wilson,
{\it Boundedness questions for Calabi-Yau threefolds},
arXiv:1706.01268.

\end{thebibliography}
\end{document}